\def\versiondate{27 Feb. 2015}
\input math10.macros
\vsize43pc
\hsize28pc

\let\nobibtex = t
\let\noarrow = t
\input eplain
\beginpackages
\usepackage{url}
   \usepackage{color}  
   \enablehyperlinks[idxexact]
\endpackages

\input Ref.macros

\checkdefinedreferencetrue
\continuousfigurenumberingtrue
\theoremcountingtrue
\sectionnumberstrue
\forwardreferencetrue
\citationgenerationtrue
\nobracketcittrue
\hyperstrue
\initialeqmacro


\input\jobname.key
\bibsty{../../texstuff/myapalike}

\def\Xsp{X}  

\def\ip#1,#2{\langle #1, #2 \rangle}
\def\Bigip#1,#2{\Big\langle #1, #2 \Big\rangle}
\def\bigip#1,#2{\big\langle #1, #2 \big\rangle}

\def\HH{{\Bbb H}}

\def\SRB{\ref b.SRB:measure/, hereinafter referred to as 
{\htmllocref{\bibcode{SRB:measure}}{SRB}}%
\def\SRB{\htmllocref{\bibcode{SRB:measure}}{SRB}}}

\ifproofmode \relax \else\head{To appear in {\it Illinois J. Math.}}
{Version of \versiondate}\fi 
\vglue20pt

\title{Hyperbolic Space Has Strong Negative Type}

\author{Russell Lyons}

\abstract{It is known that hyperbolic spaces have strict negative type, a
condition on the distances of any finite subset of points. We show that
they have strong negative type, a condition on every probability
distribution of points (with integrable distance to a fixed point). This
implies that the function of expected distances to points determines the
probability measure uniquely. It also implies that the distance covariance
test for stochastic independence,
introduced by Sz\'ekely, Rizzo and Bakirov, is consistent against all
alternatives in hyperbolic spaces.
We prove this by showing an analogue of the Cram\'er-Wold device.}

\bottomIII{Primary 
51K99,    
51M10.  
Secondary
30L05,  	
53C20.  
}
{Cram\'er-Wold, expected distances.}
{Research partially supported by NSF grant DMS-1007244 and Microsoft
Research.}

\bsection{Introduction}{s.intro}

Let $(X, d)$ be a metric space. One says that $(X, d)$ has \dfn{negative
type} if for all $n \ge 1$ and all lists of
$n$ red points $x_i$ and $n$ blue points $x'_i$ in $X$,
the sum $2\sum_{i, j} d(x_i, x'_j)$ of the distances between the $2n^2$
ordered pairs of points of opposite color is at least the sum $\sum_{i, j}
\big(d(x_i, x_j) + d(x'_i, x'_j)\big)$ of the distances between the
$2n^2$ ordered pairs of points of the same color.
It is not obvious that euclidean space has this property, but it is well
known.
By considering repetitions of $x_i$ and taking limits, we arrive at a
superficially more general property:
For all $n \ge 1$, $x_1, \ldots, x_n \in \Xsp$, and $\alpha_1, \ldots,
\alpha_n
\in \R$ with $\sum_{i=1}^n \alpha_i = 0$, we have 
$$
\sum_{i, j \le n} \alpha_i \alpha_j d(x_i, x_j) \le 0
\,.
\label e.ntdef
$$
We say that $(\Xsp, d)$ has \dfn{strict negative type} if, for every $n$
and all $n$-tuples of distinct points $x_1, \ldots, x_n$,
equality holds in \ref e.ntdef/ only when $\alpha_i = 0$ for all $i$.
Again, euclidean spaces have strict negative type.
A simple example of a metric space of non-strict negative type is $\ell^1$
on a 2-point space, i.e., $\R^2$ with the $\ell^1$-metric.

A (Borel) probability measure $\mu$ on $X$ has \dfn{finite first moment} if
$\int d(o, x) \,d\mu(x) < \infty$ for some (hence all) $o \in X$; write
$P_1(X, d)$ for the set of such probability measures.
Suppose that $\mu_1, \mu_2 \in P_1(X, d)$.
By approximating $\mu_i$ by probability measures of finite support,
we obtain a yet more general property, namely,
that when $\Xsp$ has negative type,
$$
\int d(x_1, x_2) \,d(\mu_1 - \mu_2)^2(x_1, x_2) \le 0
\,.
\label e.negtype
$$
We say that $(\Xsp, d)$ has \dfn{strong negative type} if it has negative
type and equality holds in \ref e.negtype/ only when $\mu_1 = \mu_2$. 
See \ref b.Lyons:dcov/ for an example of a (countable) metric space of
strict but not strong negative type.
The notion of strong negative type was first defined by \ref b.ZKK/.
\ref b.Lyons:dcov/
used it to show that a metric space $X$ has strong negative type
iff the theory of distance covariance holds in $X$ just as in euclidean
spaces, as introduced by \ref b.SRB:measure/.
\ref b.Lyons:dcov/ noted that
if $(\Xsp, d)$ has negative type, then
$(\Xsp, d^r)$ has strong negative type when $0 < r < 1$.

Define
$$
a_\mu(x) := \int d(x, x') \,d\mu(x')
$$
for $x \in X$ and $\mu\in P_1(X, d)$.
\ref b.Lyons:dcov/ remarked that if $(X, d)$ has negative type, then 
the map $\alpha \colon \mu \mapsto a_\mu$ is
injective on $\mu \in P_1(\Xsp)$
iff $\Xsp$ has strong negative type.
(There are also metric spaces not of negative type for which
$\alpha$ is injective.)

The concept of negative type is old, but has enjoyed a resurgence of
interest recently due to its uses in theoretical computer science, where
embeddings of metric spaces, such as graphs, play a useful role in
algorithms; see, e.g., \ref b.Naor:ICM/ and \ref b.DezaLaurent/.
A list of metric spaces of negative type appears as Theorem 3.6 of \ref
b.Meckes:pdms/; in particular, this includes
all $L^p$ spaces for $1 \le p \le 2$.
On the other hand, $\R^n$ with the $\ell^p$-metric
is not of negative type whenever $3 \le n \le \infty$
and $2 < p \le \infty$, as proved by 
\ref b.Dor/ combined with Theorem 2 of \ref b.BDK/;
see \ref b.KolLon/ for an extension to spaces that include some Orlicz
spaces, among others.
\refbmulti{Schoenberg:Ann,Schoenberg:TAMS} showed that
$X$ is of negative type iff there is a Hilbert space $H$ and a
map $\phi \colon X \to H$ such that $\all {x, x' \in X} d(x, x') =
\|\phi(x) - \phi(x')\|^2$.

That real and complex
hyperbolic spaces $\HH^n$ have negative type was shown by \ref
b.Gangolli/, Sec.~4, and was made explicit by \ref b.FarHarz/, Corollary 7.4; 
that they have strict negative type was shown by \ref b.HKM/.
(The proof of those last authors has some minor errors that are easily
corrected.)
We extend this as follows:

\proclaim Theorem.
For all $n \ge 1$, real hyperbolic space $\HH^n_\R$
of dimension $n$ has strong
negative type.

It is open whether complex hyperbolic spaces $\HH^n_\C$ have strong
negative type, which would imply our theorem.
More generally, it is open whether all Cartan-Hadamard manifolds have
strong negative type, but \ref b.HKM/ showed that those that have negative
type have strict negative type. It is known that Cartan-Hadamard surfaces
have negative type: see \ref b.CCN:Busemann/.

\bsection{Proof of the Theorem}{s.proof}

Fix $o \in \HH_\R^n$.
Let
$\sigma$ be the (infinite) Borel measure on geodesic
closed half-spaces $S \subset\HH_\R^n$
that is invariant under isometries, normalized so that
$$
\sigma\big(\{o \in S, x \notin S\}\big) = d(o, x)/2
\,;
\label e.cut
$$
see \ref b.Robertson/.
Now let $\phi(x)$ be the function $S \mapsto \I S(o) - \I S(x)$ 
in $L^2(\sigma)$. It clearly satisfies Schoenberg's condition that
$d(x, y) = \|\phi(x) - \phi(y)\|^2$. 
We call this the \dfn{Crofton embedding}, as \ref
b.Crofton/ was the first to give a formula for the distance of points in
the plane in terms of lines intersecting the segment joining them.
Thus, $\HH_\R^n$ has negative type. 
In fact, we shall not use Schoenberg's theorem, even though this half is
easy.

Instead, note that for $\mu_1, \mu_2 \in P_1(\HH_\R^n)$, we have 
$$\displaylines{
\int d(x_1, x_2) \,d(\mu_1 - \mu_2)^2(x_1, x_2) 
=
\hfill\break\cr\hfill
\int \!\! \int \big|\I S(x_1) - \I S(x_2)\big|^2
\,d(\mu_1 - \mu_2)^2(x_1, x_2) \, d\sigma(S)
\,.
\cr}
$$
Expanding the square and using the facts that 
$$
\int \I S(x) \,d\nu^2(x, y) = \nu(S) \nu(X)
$$
and
$$
\int \I S(x) \I S(y) \,d\nu^2(x, y) = \nu(S)^2
\,,
$$
we obtain that
$$
\int d(x_1, x_2) \,d(\mu_1 - \mu_2)^2(x_1, x_2) 
=
-2 \int \big(\mu_1(S) - \mu_2(S)\big)^2\, d\sigma(S)
\,.
$$
This clearly proves negative type; also, 
it is easy to prove strict negative type from this, using the fact that
every finite set has a point that is in a half-space that does not contain
any other point of the set.
In order to prove strong negative type, it clearly suffices to show that
if $\mu_1(S) = \mu_2(S)$ for $\sigma$-a.e.\ $S$ and $\mu_1, \mu_2 \in 
P_1(\HH_\R^n)$, then $\mu_1 = \mu_2$.
Consider the Klein model of $\HH_\R^n$ in which the space is the open unit
ball of $\R^n$ and in which geodesics are euclidean straight lines, whence
hyperbolic half-spaces are the intersections of euclidean half-spaces with
the unit ball.
Every probability measure in the Klein model thus is a probability measure
on $\R^n$ that happens to be carried by the unit ball.
The Cram\'er-Wold device (pp.~382--3 of \ref b.Billingsley/)
now provides the desired conclusion. 
(The usual statement of the device is that if $\mu_1$ and $\mu_2$ are
probability measures on euclidean space $\R^n$ that satisfy
$\mu_1(S) = \mu_2(S)$ for {\it all} half-spaces $S$, then $\mu_1 = \mu_2$.
Its proof extends easily to the weaker hypothesis that
$\mu_1(S) = \mu_2(S)$ for all $S$ of the form $S = \{x \in \R^n \st x \cdot
t \le \alpha\}$ for a set $B$ of pairs $(t, \alpha) \in \R^n \times\R$ with
the projection $\pi B$ of $B$ to $\R^n$ being
dense in $\R^n$ and for each $t \in \pi B$, the set $\{\alpha \st (t,
\alpha) \in B\}$ being dense in $\R$.
Alternatively, we may appeal to the fact that $\R^n$ has strong
negative type for our desired conclusion.)
\Qed

\medbreak
\noindent {\bf Acknowledgements.}\enspace 
I thank Thang Nguyen and Rich\'ard Balka for discussions.

\def\noop#1{\relax}
\input \jobname.bbl

\filbreak
\begingroup
\eightpoint\sc
\parindent=0pt\baselineskip=10pt

Department of Mathematics,
831 E. 3rd St.,
Indiana University,
Bloomington, IN 47405-7106
\emailwww{rdlyons@indiana.edu}
{http://pages.iu.edu/\string~rdlyons/}

\endgroup

\bye